\documentclass[a4,twoside,10pt]{amsart}

\usepackage{a4wide}
\usepackage{amsfonts}
\usepackage{amsmath}
\usepackage{amssymb}
\usepackage{amsthm}
\usepackage[english]{babel}
\usepackage{enumerate}
\usepackage{fontenc}
\usepackage[colorlinks=true,linkcolor=blue,citecolor=red,urlcolor=blue]{hyperref}
\usepackage[latin1]{inputenc}
\usepackage{orcidlink}
\usepackage{stackrel}
\usepackage[english]{varioref}
\usepackage[all]{xy}
\usepackage{url}

\title{Curve semistable Higgs bundles and smooth projective varieties whose canonical bundle is ample}

\author{Armando Capasso}
\address{Universit\`a degli Studi di Trieste, P.le Europa 1, Trieste (Italy), C.A.P. 34127}
\email{armando.capasso@units.it}

\thanks{A.C. is member of INdAM-GNSAGA. ORCID: 0009-0001-5463-7221 \orcidlink{0009-0001-5463-7221}}


\date{\today}

\theoremstyle{plain}
\newtheorem{theorem}{Theorem}[section]

\newtheorem{corollary}[theorem]{Corollary}
\newtheorem{lemma}[theorem]{Lemma}
\newtheorem{proposition}[theorem]{Proposition}
\newtheorem{propositions}[theorem]{Propositions}

\theoremstyle{definition}
\newtheorem{conjecture}[theorem]{Conjecture}
\newtheorem{definition}[theorem]{Definition}
\newtheorem{ex}[theorem]{Example}
\newenvironment{example}{\begin{ex}}{\hfill{$\triangle$}\end{ex}}
\newtheorem{rem}[theorem]{Remark}
\newenvironment{remark}{\begin{rem}}{\hfill{$\Diamond$}\end{rem}}
\newtheorem{remarks}[theorem]{Remarks}
\newtheorem*{Proof}{Proof}
\newenvironment{prf}{\begin{Proof}}{\hfill\text{(Q.e.d.)}\end{Proof}}

\DeclareMathOperator{\Gr}{Gr}
\DeclareMathOperator{\hGr}{\mathfrak{Gr}}
\let\hom\relax
\DeclareMathOperator{\hom}{Hom}
\DeclareMathOperator{\Id}{Id}

\DeclareMathOperator{\Pic}{Pic}
\DeclareMathOperator{\rank}{rank}
\DeclareMathOperator{\spec}{Spec}

\newcommand{\wtC}{\widetilde{C}}
\newcommand{\wvarphi}{\widetilde{\varphi}}

\newcommand{\B}{\mathbb{B}}
\newcommand{\C}{\mathbb{C}}
\newcommand{\K}{\mathbb{K}}
\newcommand{\bP}{\mathbb{P}}
\newcommand{\Q}{\mathbb{Q}}
\newcommand{\R}{\mathbb{R}}
\newcommand{\U}{\mathbb{U}}
\newcommand{\V}{\mathbb{V}}
\newcommand{\W}{\mathbb{W}}
\newcommand{\Z}{\mathbb{Z}}

\newcommand{\cE}{\mathcal{E}}
\newcommand{\cF}{\mathcal{F}}
\newcommand{\cK}{\mathcal{K}}
\newcommand{\cO}{\mathcal{O}}
\newcommand{\cQ}{\mathcal{Q}}

\newcommand{\fE}{\mathfrak{E}}
\newcommand{\fF}{\mathfrak{F}}
\newcommand{\fK}{\mathfrak{K}}
\newcommand{\fL}{\mathfrak{L}}
\newcommand{\fQ}{\mathfrak{Q}}
\newcommand{\fS}{\mathfrak{S}}

\newcommand{\rA}{\mathrm{A}}
\newcommand{\rH}{\mathrm{H}}
\newcommand{\rN}{\mathrm{N}}

\newcommand{\olC}{\overline{C}}

\newcommand{\ue}{\underline{e}}
\newcommand{\uzero}{\underline{0}}

\newcommand{\wtvarphi}{\widetilde{\varphi}}

\begin{document}

\begin{abstract}
Considering the so-called Simpson system on smooth projective varieties, defined over an algebraically closed field of characteristic $0$, whose canonical bundle is ample, we provide a new proof of the stability of this Higgs bundle, which in turn yields another proof of the Miyaoka-Yau inequality. When equality holds, we prove that the discriminant class of the Simpson system is numerically trivial and that this Higgs bundle is curve semistable. This result follows from a study of the relations between weak ampleness, numerically effectiveness for Higgs bundles that ``feel'' the Higgs field and (semi)stability. Moreover, we obtain another proof of the algebraic hyperbolicity of these varieties which provides a lower bound on a positive real constant related to this property; to best of our knowledge, this is the first and only result of this type. We conjecture that a similar lower bound holds for smooth projective varieties whose canonical bundle is ample.
\end{abstract}

\maketitle

\markboth{Curve semistable Higgs bundles and smooth projective varieties whose canonical bundle is ample}{Armando Capasso}

\section*{Introduction}

\noindent Let $X$ be a smooth, complex, projective variety of dimension $n\geq2$. Recall that for any ${k\in\{1,\dotsc,n\}}$, $c_k(X)=(-1)^kc_k\left(\Omega_X^1\right)$, where $\Omega^1_X$ is the \emph{cotangent bundle of} $X$; under the assumption that $K_X$ (the \emph{canonical bundle}) is \emph{ample}, Yau proved the inequality
\begin{equation}\label{eqMY}
MY(X)\stackrel{def.}{=}(-1)^n\int_X\left[c_2(X)-\frac{n}{2(n+1)}c_1(X)^2\right]\cdot c_1(X)^{n-2}\geq0,
\end{equation}
and the equality holds if and only if $X$ is covered by the complex ball $\B^n$ holomorphically \cite[Theorem 4 and Remark (iii)]{Y:S-T}. Moreover, assuming $n=2$ and that $X$ saturates \eqref{eqMY}, Miyaoka has proved that $K_X$ and $\Omega^1_X$ are both ample \cite[Corollary at page 294]{M:Y:3}. On the other hand, under the hypothesis that $K_X$ is big and numerically effective, Miyaoka proved the same inequality for $n=2$ \cite[Theorem 4]{M:Y:2}.\medskip

\noindent Following another way, Simpson introduced an opportune \emph{Higgs bundle} $\fS=(S,\varphi)$ (see \eqref{eqSS}) over $X$, and he proved the stability of $\fS$ under suitable hypotheses on $X$ (see \cite[Proposition 9.9]{S:CT:1}). By stability of $\fS$ with respect to a polarization $L$ of $X$ (\cite[Proposition 9.8]{S:CT:1}), one has:
\begin{displaymath}
(-1)^n\int_X\left[c_2(X)-\frac{n}{2(n+1)}c_1(X)^2\right]\cdot L^{n-2}=\int_X\frac{1}{2(n+1)}c_2\left(S\otimes S^{\vee}\right)\cdot L^{n-2}\geq0;
\end{displaymath}
if $L$ is the canonical bundle, one recovers the inequality \eqref{eqMY}. For simplicity, the characteristic class $\displaystyle\Delta(S)=\frac{1}{2(n+1)}c_2\left(S\otimes S^{\vee}\right)\in H^4(X,\Q)$ is called \emph{discriminant of $S$}.
\begin{remarks}
To be clear:
\begin{enumerate}[a)]
\item \cite[Proposition 9.8]{S:CT:1} assumes \emph{a priori} that $\fS$ is stable with respect to a polarization $L$ of $X$, without any assumption on the dimension of $X$;
\item instead, among the various hypotheses, in \cite[Proposition 9.9]{S:CT:1} one assumes that $X$ has dimension $2$. \hfill{$\Diamond$}
\end{enumerate}
\end{remarks}
\noindent These results have been extended by Greb, Kebekus, Peternell and Taji to suitable singular, complex, projective varieties (\cite[Example 5.3 and Corollary 7.2]{GKPT}). In particular, from all this one obtains a generalisation of \eqref{eqMY} to this general setting (\cite[Theorem 1.1]{GKPT}).\medskip
%
%\noindent Simpson, using an opportune \emph{Higgs bundle} $\fS=(S,\varphi)$ (see \eqref{eqSS}), generalized Yau's result to compact K\"ahler manifolds $(M,\omega)$ \cite[Propositions 9.8 and 9.9]{S:CT:1}. This is accomplished proving the $\omega$-\emph{stability} of that Higgs bundle which implies
%\begin{displaymath}
%GY(M)=\int_M\left[c_2(M)-\frac{n}{2(n+1)}c_1(M)^2\right]\cdot\omega^{n-2}=\int_M\frac{1}{2(n+1)}c_2\left(S\otimes S^{\vee}\right)\cdot\omega^{n-2}\geq0.
%\end{displaymath}
%The characteristic class $\displaystyle\Delta(S)=\frac{1}{2(n+1)}c_2\left(S\otimes S^{\vee}\right)\in H^2(M,\Q)$ is called \emph{discriminant of $S$}.\medskip

\noindent From now on, let $X$ be defined over an algebraically closed field of characteristic $0$. Under the hypothesis that $n\geq2$, $K_X$ is big and numerically effective, Miyaoka proved an inequality similar to \eqref{eqMY} \cite[Proposition 7.1]{M:Y:1}.\medskip

\noindent Assuming the ampleness of $K_X$, we provide another proof the stability of $\fS$ with respect to $K_X$ (Theorem \ref{th3.1}) from which another proof follows of the inequality \eqref{eqMY} (Theorem \ref{th3.3}). In the extremal case, \emph{i.e.} when $\displaystyle\int_X\Delta(S)\cdot K_X^{n-2}=0$, we prove that $\fS$ is \emph{curve semistable} (Definition \ref{def1.1} and Corollary \ref{cor3.1}.\ref{cor3.1.a}) hence it is semistable with respect to any polarization of $X$ (Corollary \ref{cor3.1}.\ref{cor3.1.c}) and the equivalence $\displaystyle c_2(X)-\frac{n}{2(n+1)}c_1(X)^2\equiv_{num}0$ holds (Corollary \ref{cor3.1}.\ref{cor3.1.b}).\medskip

\noindent The study of curve semistable Higgs bundles, in the complex setting, started by Bruzzo and Hern\'andez Ruip\'erez in \cite{B:HR} and it had been refined by Bruzzo and Gran\~na Otero in \cite{B:GO:1,B:GO:2}. In \cite{A:C}, we proved that the study of curve semistability for Higgs bundles over smooth projective varieties, defined over an algebraically closed field of characteristic $0$, may be reduced to complex setting by suitable base changes (Corollary 2.10 \emph{ibidem}).\medskip

\noindent On the other hand, in \cite{B:GO:1,B:GO:2,B:M:R} the authors defined notions of \emph{ampleness} and \emph{numerical effectiveness} for Higgs bundles (Definitions \ref{def2.1}, \ref{def2.3} and \ref{def2.2}), using a generalization of the Grassmann bundles of a vector bundle, called the \emph{Higgs Grassmann schemes}, introduced in \cite{B:HR}.
\begin{remark}\label{rem0.1}
The notions of ampleness for Higgs bundle recalled above are not equivalent (\cite[Remark 4]{B:M:R}). In order to avoid confusion, the ampleness condition for Higgs bundles introduced in \cite{B:GO:1,B:GO:2} here is called \emph{strong Higgs ampleness} (cfr. Definition \ref{def2.3}), whereas the ampleness condition for Higgs bundles introduced in \cite{B:M:R} here is called \emph{weak Higgs ampleness} (cfr. Definition \ref{def2.2}).
\end{remark}
\noindent The underlying idea is to formulate notions of positivity that ``feel'' the Higgs field, and reduce to the classical ones when the Higgs field vanishes. This holds for the \emph{Higgs numerical effectiveness} and \emph{weak Higgs ampleness} (Remark \ref{rem2.1}.\ref{rem2.1.a}) but not for \emph{strong Higgs ampleness} (Remark \ref{rem2.3}). In this paper we present an explicit example of a \emph{weakly H-ample} Higgs bundle that is not \emph{strongly H-ample} (see Proposition \ref{prop3.1}, cfr. also \cite[Remark 4]{B:M:R})\medskip

\noindent From Section \ref{SSMinVarGenType} onward, assuming in addition that $X$ saturates \eqref{eqMY}, by curve semistability of Simpson system $\fS$ implies a further result on $\fS$ (Theorem \ref{th4.1}, cfr. \cite[Proposition 5 at page 293]{M:Y:3}) which imposes certain restrictions on the curves in $X$ (Corollary \ref{cor4.1}). Since $K_X$ is ample, the weak H-ampleness of $\fS$ also follows (Lemma \ref{lem3.1}) by applying \cite[Corollary 3.4]{B:M:R}. Moreover, since $\left(\Omega^1_X,0\right)$ is a quotient Higgs bundle of $\fS$, it is ample (Corollary \ref{cor3.2}).\medskip

\noindent Moreover, in the complex setting, the stability of $\fS$ allows us to give a new proof of the fact that if $X$ saturates \eqref{eqMY}, then $X$ is uniformized by the complex ball $\B^n$ (Theorem \ref{th5.1}); even though there are general results in suitable singular, complex, projective settings (\cite[Theorems 1.2 and 1.3]{GKPT}). On the other hand, the curve semistability of $\fS$ allows us to give a new proof of the \emph{algebraic hyperbolicity} of $X$ (Definition \ref{def5.1} and Theorem \ref{th3.5}). The advantage of this new proof is that it yields a lower bound for a positive real constant depending on $X$ and $K_X$: to best of our knowledge, this is the first and only result of this type. We conjecture that a similar result holds for any smooth, complex, projective variety whose canonical bundle is ample (Conjecture \ref{conj}).\medskip

\noindent \textbf{Notations and conventions.} $\K$ is an algebraically closed field of characteristic $0$, unless otherwise indicated. By a projective variety $X$ we mean a projective, integral scheme of finite type over $\K$ of dimension $n\geq1$. If $n\in\{1,2\}$ we shall write projective curve or projective surface, respectively. Whenever we consider a morphism $f\colon C\to X$, it is understood that $C$ is an irreducible smooth projective curve.

\section{Curve semistable Higgs bundles}

\noindent Let $X$ be a smooth scheme over $\K$, let $\Omega^1_X$ be the cotangent bundle of $X$.
\begin{definition}
A \emph{Higgs sheaf} $\fE$ on $X$ is a pair $(\cE,\varphi)$ where $\cE$ is an $\cO_X$-coherent sheaf equipped with a morphism $\varphi\colon\cE\to\cE\otimes\Omega^1_X$ called \emph{Higgs field} such that the composition
\begin{displaymath}
\varphi\wedge\varphi\colon\cE\xrightarrow{\varphi}\cE\otimes\Omega^1_X\xrightarrow{\varphi\otimes\Id}\cE\otimes\Omega^1_X\otimes\Omega^1_X\to\cE\otimes\Omega^2_X
\end{displaymath}
vanishes. A Higgs subsheaf of $\fE$ is a $\varphi$-\emph{invariant subsheaf} $\cF$ of $\cE$, \emph{i.e.} $\varphi(\cF)\subseteq\cF\otimes\Omega_X^1$. A \emph{quotient Higgs sheaf} of $\fE$ is a quotient sheaf of $\cE$ such that the corresponding kernel is $\varphi$-invariant. A \emph{Higgs bundle} is a Higgs sheaf whose underlying coherent sheaf is locally free.
\end{definition}
\noindent Let $H$ be a polarization of $X$ and let $\fE=(\cE,\varphi)$ be a torsion-free Higgs sheaf on $X$, if not otherwise indicated.
One defines the \emph{slope} of $\fE$ as $\displaystyle\mu(\fE)=\frac{1}{\rank(\cE)}\int_Xc_1(\det(\cE))\cdot H^{n-1}\in\Q$.
\begin{definition}
$\fE$ is $H$-(\emph{semi})\emph{stable} if $\mu(\fF)\stackrel[(\leq)]{\textstyle<}{}\mu(\fE)$ for every Higgs subsheaf $\fF$ of $\fE$ with ${0<\rank(\fF)<\rank(\fE)}$, equivalently, if $\mu(\fE)\stackrel[(\leq)]{\textstyle<}{}\mu(\fQ)$ for every torsion-free quotient Higgs sheaf $\fQ$ of $\fE$ with $0<\rank(\fQ)<\rank(\fE)$. In the other eventuality, $\fE$ is \emph{unstable}.
\end{definition}
\noindent We recall that every torsion-free Higgs sheaves $\fE$ has a unique \emph{Harder-Narasimhan filtration} (\emph{HN-filtration}, for short)
\begin{equation}\label{eqHNfiltr}
{0}=\fE_0\subsetneqq\fE_1\subsetneqq\dotsc\subsetneqq\fE_{m-1}\subsetneqq\fE_m=\fE
\end{equation}
such that the successive quotient Higgs sheaves $\fQ_i=\fE_i/\fE_{i-1}$ are $H$-stable and their slopes satisfy the condition
\begin{displaymath}
\mu_{\max}(\fE)\stackrel{def.}{=}\mu(\fQ_1)>\dotsc>\mu(\fQ_m)\stackrel{def.}{=}\mu_{\min}(\fE).
\end{displaymath}
For every quotient Higgs sheaf $\fQ$ of $\fE$ one has $\mu(\fQ)\geq\mu_{\min}(\fE)$.
\begin{definition}\label{def1.1}
$\fE$ is \emph{curve semistable} if $f^{\ast}\fE$ is semistable for any $f\colon C\to X$.
\end{definition}
\noindent For simplicity, we shall skip any reference to the fixed polarization $H$ of $X$ if there is no confusion.\medskip

\noindent We advice that, by \cite[Theorem 1.9 and Section 11]{M:YI}, we can apply the usual theory of Chern classes for locally free sheaves to coherent sheaves on $X$. Thus we consider the characteristic class
\begin{displaymath}
\Delta(\cE)=\frac{1}{2r}c_2\left(\cE\otimes\cE^{\vee}\right)=c_2(\cE)-\frac{r-1}{2r}c_1(\cE)^2\in\rA^2(X)\otimes_{\Z}\Q=\rA^2(X)_{\Q},
\end{displaymath}
called the \emph{discriminant} of $\cE$. Here $r=\rank(\cE)$, $\rA_k(X)$ is the Abelian group of $k$-cycles on $X$ modulo rational equivalence and $\rA^k(X)=\rA_{n-k}(X)$.\medskip

\noindent Henceforth, let $X$ be a smooth projective variety.
\begin{theorem}[{\cite[Theorem 7]{L:A}}]\label{th1.1}
Let $\fE=(\cE,\varphi)$ be a semistable Higgs sheaf on $X$ with respect to a polarization $H$. Then
\begin{displaymath}
\int_X\Delta(\cE)\cdot H^{n-2}\geq0.
\end{displaymath}
\end{theorem}
\begin{remark}
If one assumes $\K=\C$, the previous theorem is \cite[Proposition 3.4]{S:CT:1}.
\end{remark}
\noindent About the ``extremal'' case, the following theorem holds.
\begin{theorem}[{see \cite[Theorems 1.2 and 1.3]{B:HR} and \cite[Proposition 3.2]{L:LG}}]\label{th1.2}
Let $\fE=(E,\varphi)$ be a Higgs bundle over $X$.
\begin{enumerate}[a)]
\item\label{th1.2.a} If $\fE$ is semistable with respect to some polarization $H$ and $\displaystyle\int_X\Delta(E)\cdot H^{n-2}=0$, then $\fE$ is curve semistable.
\item\label{th1.2.b} If $\fE$ is curve semistable, then $\fE$ is semistable with respect to any polarization $H$.
\end{enumerate} 
\end{theorem}
\newpage
\begin{remarks}\label{rem1.1}
\,\begin{enumerate}[a)]
\item\label{rem1.1.a} In statement \ref{th1.2.a} of the previous Theorem, one may write $\Delta(E)\equiv_{num}0$ instead of $\displaystyle{\int_X\Delta(E)\cdot H^{n-2}=0}$. Indeed, this last equality is equivalent to $\displaystyle\int_Xc_2\left(E\otimes E^{\vee}\right)\cdot H^{n-2}=0$. By hypotheses, this allows us to apply \cite[Corollary 6, see also the note after Theorem 11]{L:A}, also because $c_1\left(E\otimes E^{\vee}\right)=0$, so that $c_k\left(E\otimes E^{\vee}\right)\equiv_{num}0$ for all $k>0$, hence $\Delta(E)=2rc_2\left(E\otimes E^{\vee}\right)\equiv_{num}0$.
\item Statement \ref{th1.2.b} of the previous theorem follows from the ``easy'' direction of \cite[Theorem 6.1]{M:R}. \hfill{$\Diamond$}%It is an open problem either to prove or to disprove the vanishing of the discriminant class of a curve semistable Higgs bundle.
\end{enumerate}
\end{remarks}

\section{Positivity conditions for Higgs bundles}

\subsection{Higgs Garssmann schemes}

\noindent Let $\fE=(E,\varphi)$ be a rank $r\geq2$ Higgs bundle over a smooth projective variety $X$, unless otherwise indicated. Let $s\in\{1,\dotsc,r-1\}$, and let $p_s\colon\Gr_s(E)\to X$ be the \emph{Grassmann bundle} parametrizing rank $s$ locally free quotients of $E$ (see \cite{F:W}). Consider the short exact sequence of vector bundles over $\Gr_s(E)$ 
\begin{displaymath}
\xymatrix{
0\ar[r] & S_{r-s,E}\ar[r]^(.55){\eta} & p_s^{\ast}E\ar[r]^(.45){\zeta} & Q_{s,E}\ar[r] & 0
},
\end{displaymath}
where $S_{r-s,E}$ is the \emph{universal rank $r-s$ subbundle} and $Q_{s,E}$ is the \emph{universal rank $s$ quotient bundle of} $p_s^{\ast}E$, respectively. One defines the closed subschemes $\hGr_s(\fE)\subseteq\Gr_s(E)$ (the $s$-\emph{th Higgs Garssmann schemes of} $\fE$) as the zero loci of the composite morphisms
\begin{displaymath}
(\zeta\otimes\Id)\circ p_s^{\ast}\varphi\circ\eta\colon S_{r-s,E}\to Q_{s,E}\otimes p_s^{\ast}\Omega_X^1.
\end{displaymath}
The restriction of the previous sequence to $\hGr_s(\fE)$ yields a universal short exact sequence
\begin{displaymath}
\xymatrix{
0\ar[r] & \mathfrak{S}_{r-s,\fE}\ar[r]^(.55){\eta} & \rho_s^{\ast}\fE\ar[r]^(.45){\zeta} & \fQ_{s,\fE}\ar[r] & 0,
}
\end{displaymath}
where $\fQ_{s,\fE}=Q_{s,E|\hGr_s(\fE)}$ is equipped with the quotient Higgs field induced by $\rho_s^{\ast}\varphi$, and ${\rho_s=p_{s|\hGr_s(\fE)}}$. The scheme $\hGr_s(\fE)$ enjoys the usual universal property: for a morphism of varieties $f\colon Y\to X$, the morphism $g\colon Y \to \Gr_s(E)$ given by a rank $s$ quotient $Q$ of $f^{\ast} E$ factors through $\hGr_s(\fE)$ if and only if $\varphi$ induces a Higgs field on $Q$.\medskip

\noindent In order to explain a motivation behind the definition of \emph{weakly H-ample Higgs bundle}, we prove the following.
\begin{proposition}[{cfr. \cite[Proposition 3.7]{B:HR}}]\label{prop2.3}
Let $\fE=(E=L_1\oplus L_2,\varphi)$ be a rank $2$ Higgs bundle over $(X,H)$ where $H$ is a polarization of $X$, such that ${\varphi(L_1)\subseteq L_2\otimes\Omega^1_X}$ and $\varphi(L_2)=0$. Then
\begin{enumerate}[a)]
\item $\displaystyle\varphi=\begin{pmatrix}
0 & 0\\
c & 0
\end{pmatrix}\neq0$ only if $\deg\left(L_1\right)\leq\deg\left(L_2\right)+\mu_{\max}\left(\Omega^1_X\right)$, where $c\in\hom_{\cO_X}\left(L_1,L_2\otimes\Omega^1_X\right)$ is a monomorphism;
\item $\hGr_1(\fE)$ coincides with $\Gr_1\left(L_1\right)$.
\end{enumerate}
\end{proposition}
\begin{prf}
For our aims we assume $\varphi\neq0$. Let
\begin{displaymath}
\varphi=\begin{pmatrix}
a & b\\
c & d
\end{pmatrix}
\end{displaymath}
where
\begin{gather*}
a\in\hom_{\cO_X}\left(L_1,L_1\otimes\Omega^1_X\right),b\in\hom_{\cO_X}\left(L_2,L_1\otimes\Omega^1_X\right),\\
c\in\hom_{\cO_X}\left(L_1,L_2\otimes\Omega^1_X\right),d\in\hom_{\cO_X}\left(L_2,L_2\otimes\Omega^1_X\right);
\end{gather*}
by hypotheses, one has $a=0$, $b=0$, $d=0$. Trivially $\left(L_2,\varphi_{\vert L_2}=0\right)$ is a Higgs subbundle of $\fE$ and thence $\fL_1=\left(L_1,\wtvarphi\right)$ is a quotient Higgs bundle of $\fE$. %Let $\xi$ be the generic point of $X$. If $\ker(c)\left(\xi\right)\neq0$ then by Lemma \ref{lemHS.4} and Proposition \ref{propDET4.1}.c $\ker(c)=L_1$ which is the case that we are excluding. Thus $\ker(\varphi)=L_2$, because $L_1$ is torsion-free and the condition $\ker(c)\left(\xi\right)=0$ implies $\ker(c)=\uzero_X$. Moreover, the following diagram
%\begin{displaymath}
%\xymatrix{
%0\ar[r] & L_2\ar[d]_0\ar@{^{(}->}[r]^i & E\ar[d]_{\varphi}\ar@{->>}[r]^p & L_1\ar[d]_{\wvarphi}\ar[r] & 0\\
%0\ar[r] & 0\ar[r] & E\otimes\Omega^1_X\ar[r]_{p\otimes Id} & L_1\otimes\Omega^1_X\ar[r] & 0
%}
%\end{displaymath}
%is commutative by definition. For all $x\in X$ let $s_x\in L_{1,x}$. It is ${p_x^{-1}(s_x)=\{(s_x,t_x)\in E_x\mid t_x\in L_{2,x}\}}$ and from this $(p\otimes\Id)_x(\varphi_x(s_x,t_x))=(p\otimes\Id)_x(0,c_x(s_x))=0$ \emph{i.e.} $\wvarphi=0$.
\smallskip
%
%\noindent Let $\psi\colon E\otimes\left(\Omega^1_X\right)^{\vee}\to E\otimes\Omega^1_X\otimes\left(\Omega^1_X\right)^{\vee}\xrightarrow{\Id_E\otimes\Tr}E$, by definition $\coker(\psi)=E/\Im(\psi)$ and we compute it as it follows
%\begin{displaymath}
%\psi((s,t)\otimes t)=\left(\Id\otimes\Tr\right)((0,c(s))\otimes t)\in L_2
%\end{displaymath}
%\emph{i.e.} $\Im(\psi)\subseteq L_2$ and by \cite[Proposition II.5.7]{H:RC:1} this is a coherent subsheaf of $L_2$, moreover, by hypothesis this is also a Higgs subsheaf of $\fE$. ...\smallskip

\noindent Let $\fK=\left(\cK,\varphi_{\vert\cK}\right)$ be a Higgs subsheaf of $\fE$ such that the corresponding quotient Higgs sheaf ${\fQ_0=\left(\cQ_0,\wvarphi_0\right)}$ is locally free and has rank $1$. Thus $\cK$ is a locally free sheaf, because it is a kernel of an epimorphism of locally free sheaves on a smooth projective variety and has rank $1$ (see for example \cite[Propositions 1.1 and 1.9]{H:RC:2}). By hypothesis ${\cK=\left(\cK\cap L_1\right)\oplus\left(\cK\cap L_2\right)\equiv\cK_1\oplus\cK_2}$, hence $\varphi\left(\cK_2\right)=0$ and $\varphi\left(\cK_1\right)\subseteq\cK_2\otimes\Omega^1_X$.\smallskip %, where $\cK_1$ and $\cK_2$ are the underlying coherent sheaves to $\fK_1$ and $\fK_2$, respectivelySince $\cK_i$'s are kernel of a morphism from a locally free sheaf on a smooth projective curve, we infer their locally freeness, where $i\in\{1,2\}$.

\noindent Let $\eta$ be the generic point of $X$ and let $\_(\eta)$ be the generic fibre of the sheaf $\_$. One has two possibilities:
\begin{enumerate}[1)]
\item $\cK_1(\eta)=\left\{\uzero\right\}$ \emph{i.e.} $\cK_1$ is a torsion sheaf, hence $\cK_1=\uzero_X$ because $L_1$ is torsion-free. It follows that $\cK_2=\cK=L_2$.
\item $\cK_1(\eta)\neq\left\{\uzero\right\}$ then $\cK_2(\eta)=\left\{\uzero\right\}$, \emph{i.e.} $\cK_2$ is a torsion sheaf. By the same reasoning $\cK_1=L_1$ and $\cK_2=\uzero_X$, this forces $\varphi$ to be zero and this is impossible because we are supposing the contrary.
\end{enumerate}
From all this, $\fL_1$ is the only quotient Higgs bundle of $\fE$ and the claim follows.
\end{prf}

\subsection{H-nef Higgs bundles}

\noindent The \emph{Higgs nefness condition} for $\fE$ amounts to saying that suitable line bundles are \emph{nef} in the usual sense (\cite[Proposition A.4]{B:GO:2}). Following \cite[Definition 1.4.1]{L:RK}, we recall that a line bundle $L$ over a proper scheme $Z$ is \emph{numerically effective} (\emph{nef}, for short) if $\displaystyle\int_Cc_1(L)\geq0$ for every irreducible curve $C\subseteq Z$.
\begin{definition}[{see \cite[Definition 2.3]{B:GO:1}}]\label{def2.1}
A Higgs bundle $\fE=(E,\varphi)$ of rank one is said to be \emph{Higgs numerically effective} (\emph{H-nef}, for short) if $E$ is nef in the usual sense. If $\rank(\fE)=r\geq2$, one defines H-nefness inductively by requiring that
\begin{enumerate}[a)]
\item all Higgs bundles $\fQ_{s,\fE}$ are H-nef for all $s\in\{1,\dotsc,r-1\}$, and
\item the determinant line bundle $\det(E)$ is nef.
\end{enumerate}
\end{definition}
\noindent Higgs bundles over $X$ which are H-nef satisfy properties analogous to those of nef vector bundles; here we list some of them.
\begin{propositions}\label{prop2.2}
\,\begin{enumerate}[a)]
\item\label{prop2.2.a} Let $f\colon Y\to X$ be a morphism of smooth projective varieties, and let $\fE$ be a Higgs bundle over $X$. If $\fE$ is H-nef then $f^{\ast}\fE$ is H-nef (\cite[Proposition 2.6.(ii)]{B:GO:1}). If $f$ is also surjective and $f^{\ast}\fE$ is H-nef then $\fE$ is H-nef (\cite[Lemma 3.4]{B:B:G}).
\item\label{prop2.2.b} Let $\fE$ be H-nef then every quotient Higgs bundle of $\fE$ is H-nef (\cite[Lemma 3.5]{B:B:G}).
\end{enumerate}
\end{propositions}
\begin{theorem}[{\cite[Theorem 3.4]{B:C:GO}}]\label{th2.2}
Let $\fE$ be a Higgs bundle over $X$. Then $\fE$ is H-nef if and only if $\mu_{\min}\left(f^{\ast}\fE\right)\geq0$ for all morphisms $f\colon C\to X$, where $\mu_{\min}\left(f^{\ast}\fE\right)$ is defined via the relevant HN-filtration \eqref{eqHNfiltr}.
\end{theorem}

\subsection{Weakly and strongly H-ample Higgs bundles}

\noindent Mimicking Definition \ref{def2.1}, one has the following one.
\begin{definition}[{cfr. \cite[Definition 2.3]{B:GO:1}}]\label{def2.3}
A Higgs bundle $\fE=(E,\varphi)$ of rank one is said to be \emph{strongly Higgs ample} (\emph{strongly H-ample}, for short) if $E$ is ample in the usual sense. If ${\rank(\fE)=r\geq2}$, one defines strongly H-ampleness inductively by requiring that
\begin{enumerate}[a)]
\item all Higgs bundles $\fQ_{s,\fE}$ are strongly H-ample for all $s\in\{1,\dotsc,r-1\}$, and
\item the determinant line bundle $\det(E)$ is ample.
\end{enumerate}
\end{definition}
\noindent As stated in Remark \ref{rem0.1}, in the previous definition we have added the term ``strongly'' with respect to the original definition, because there is another notion of ampleness for Higgs bundles. The following remark explains a motivation behind this choice; another motivation is explained in Remark \ref{rem2.1}.\ref{rem2.1.a}.
\begin{remark}\label{rem2.3}
Let $\fE=(E,\varphi)$ be a rank $r\geq3$ strongly H-ample Higgs bundle over $X$ then $\varphi\neq0$, because $Q_{s,E}$ is never ample for $s\in\{2,\dotsc,r-1\}$ (cfr. \cite[Example 6.1.6]{L:RK}).\smallskip

\noindent To be clear, let $x\in X$ be a closed point. $\Gr_s(E)_x=\Gr_s\left(E_x\right)$ is isomorphic to the \emph{Grassmannian} $\Gr_s\left(\V^{\vee}\right)$ \emph{of} $s$-\emph{dimensional subspaces of the dual of a} $r$-\emph{dimensional vector space} $\V$, and this last one is a closed subscheme of $\Gr_s(E)$, up to closed embeddings. Let $\left\{\ue_1,\dotsc,\ue_r\right\}$ be a basis of $\V$; let $\W\subsetneqq\U\subseteq\V^{\vee}$ be subspaces generated by $\left\{\ue_1^{\vee},\dotsc,\ue_{s-1}^{\vee}\right\}$ and $\left\{\ue_1^{\vee},\dotsc,\ue_{s+1}^{\vee}\right\}$, respectively. Let $\gamma\colon\bP_{\K}^1\to\Gr_s\left(\V^{\vee}\right)$ be a map of schemes such that on closed points $[a:b]\in\bP_{\K}^1$ it results $\gamma([a:b])=\left[\left\langle\ue_1^{\vee},\dotsc,\ue_{s-1}^{\vee},a\ue_s^{\vee}+b\ue_{s+1}^{\vee}\right\rangle\right]\in\Gr_s\left(\V^{\vee}\right)$. By \cite[Theorem II.7.1.(2)]{H:RC:1}, $\gamma$ is a well-defined finite morphism of schemes and it describes a rational curve into $\Gr_s\left(\V^{\vee}\right)$. Let $Q_s$ be the universal quotient bundle over $\Gr_s\left(\V^{\vee}\right)$. $\gamma^{\ast}Q_s$ is a rank $s$ vector bundle over $\bP^1_{\K}$ hence it is direct sum $\cO_{\bP_{\K}^1}\left(d_1\right)\oplus\dotsc\oplus\cO_{\bP_{\K}^1}\left(d_s\right)$ of $s$ line bundles, unique up to permutations (\cite[Th\'eor\`eme 2.1]{G:A}). Since first $s-1$ sections are constants then they generate trivial line bundles over $\bP_{\K}^1$, that is $d_1=\dotsc=d_{s-1}=0$; instead, the last section generates $\cO_{\bP_{\K}^1}(1)$, that is $d_s=1$. By all this, $\gamma^{\ast}Q_s$ has trivial quotient bundles and this is not possible for the finite pull-backs of ample vector bundles.
\end{remark}
\noindent To address the discrepancy mentioned in the previous remark, there is the following definition.
\begin{definition}[{cfr. \cite[Definition 2.5]{B:M:R}}]\label{def2.2}
A rank $r$ Higgs bundle $\fE=(E,\varphi)$ is said to be \emph{weakly Higgs ample} (\emph{weakly H-ample}, for short) if:
\begin{enumerate}[a)]
\item $\fE$ is H-nef,
\item all line bundles $\det\left(\fQ_{s,\fE}\right)$ are ample for all $s\in\{1,\dotsc,r-1\}$, and
\item the determinant line bundle $\det(E)$ is ample.
\end{enumerate}
\end{definition}
\noindent In the previous Definitions \ref{def2.1}, \ref{def2.3} and \ref{def2.2}, the condition on the determinant cannot be omitted as \cite[Example 2.5]{B:GO:1} and the following example show.
\begin{example}[{cfr. \cite[Example 2.5]{B:GO:1}}]\label{exPCHB.1}
Let $X$ be a smooth projective curve of genus $g\geq2$, and let $\fE=(E=L_1\oplus L_2,\varphi)$ be a rank $2$ Higgs bundle over $X$ such that $\varphi(L_1)\subseteq L_2\otimes\Omega^1_X$ and $\varphi(L_2)=0$. By Proposition \ref{prop2.3}, $\fE$ has only a quotient Higgs bundle which is $\left(L_1,\widetilde{\varphi}\right)$. If we take $\deg(L_1)=1$ and $\deg(L_2)=-2$, then $\fE$ is a Higgs bundle such that $\deg(E)=-1$ and $\fQ_{1,\fE}$ is an ample line bundle.
\end{example}
\begin{remarks}\label{rem2.1}
\,
\begin{enumerate}[a)]
\item\label{rem2.1.a} Let $\fE=(E,\varphi)$ be a Higgs bundle over $X$. When $\varphi=0$, $\fE$ is weakly H-ample (H-nef, respectively) if and only if $E$ is ample (nef, respectively) as ordinary vector bundle (see also \cite[Remark 3]{B:M:R}).
\item A strongly H-ample Higgs bundle is weakly H-ample, but the converse does not hold. We refer the reader to Proposition \ref{prop3.1} for explicit examples of weakly but not strongly H-ample Higgs bundle with a Higgs field different from zero. \hfill{$\Diamond$}
\end{enumerate}
\end{remarks}
\noindent Also Higgs bundles over $X$ which are strongly or weakly H-ample satisfy properties analogous to those of ample vector bundles; here we list some of them.
\begin{propositions}\label{prop2.1}
\,
\begin{enumerate}[a)]
\item Let $f\colon Y\to X$ be a finite morphism of smooth projective varieties. If $\fE$ is strongly H-ample then $f^{\ast}\fE$ is strongly H-ample. Moreover, if $f$ is also surjective and $f^{\ast}\fE$ is strongly H-ample then $\fE$ is strongly H-ample (\cite[Proposition 2.6.(1)]{B:C:GO}).
\item Let $\fE$ be strongly H-ample then every quotient Higgs bundle of $\fE$ is strongly H-ample (\cite[Proposition 2.6.(2)]{B:C:GO}).
\item\label{prop2.1.c} Let $f\colon Y\to X$ be a finite morphism of smooth projective varieties. If $\fE$ is weakly H-ample then $f^{\ast}\fE$ is weakly H-ample. Moreover, if $f$ is also surjective and $f^{\ast}\fE$ is weakly H-ample then $\fE$ is weakly H-ample (\cite[Proposition 3.1]{B:M:R}).
\item\label{prop2.1.d} Let $\fE$ be weakly H-ample then every quotient Higgs bundle of $\fE$ is weakly H-ample (\cite[Proposition 3.2]{B:M:R}).
%\item\label{prop2.1.c} Let $\fE=(E,\varphi)$ be curve semistable and let $c_1(E)$ be ample. Then $\fE$ is H-ample (\cite[Proposition 3.11]{B:C:GO}).
\end{enumerate}
\end{propositions}
\begin{theorem}[{\cite[Theorem 3.3]{B:M:R}}]\label{th2.1}
Let $\fE=(E,\varphi)$ be a Higgs bundle over $X$. Fix an ample class $h\in\rN^1(X)$. Then $\fE$ is weakly H-ample if and only if
\begin{enumerate}[a)]
\item the line bundle $\det(E)$ is ample;
\item there exists $\delta\in\R_{>0}$ such that for every finite morphism $f\colon C\to X$, the inequality
\begin{displaymath}
\mu_{\min}\left(f^{\ast}\fE\right)\geq\delta\int_Cf^{\ast}h
\end{displaymath}
holds, where $\mu_{\min}\left(f^{\ast}\fE\right)$ is defined via the relevant HN-filtration \eqref{eqHNfiltr}.
\end{enumerate}
\end{theorem}

\section{On the Simpson system on smooth projective varieties\\ with ample canonical bundle}\label{SSMinVarGenType}

\noindent From now on, let $X$ be a smooth projective variety of dimension $n\geq2$ with $K_X$ (the \emph{canonical bundle}) ample, unless otherwise indicated. Let $\fS=(S,\varphi)$ be the \emph{Simpson system}, where ${S=\Omega^1_X\oplus\cO_X}$ and
\begin{equation}\label{eqSS}
\varphi=\begin{pmatrix}
0 & 0\\
\Id & 0
\end{pmatrix},\,\Id\in\hom\left(\Omega^1_X,\Omega^1_X\right).
\end{equation}
Actually the Simpson system is $K_X$-stable.
\begin{theorem}
The cotangent bundle $\Omega^1_X$ is $K_X$-semistable.
\end{theorem}
\begin{prf}
By \cite[Lemma 2.2.1]{AC:PhD}, there exist an algebraically closed subfield $\K_0$ of $\K$ and a smooth projective variety $X_{\K_0}$ over $\K_0$ such that $\K_0$ is isomorphic to a subfield of $\C$ and the following diagram is Cartesian
\begin{displaymath}
\xymatrix{
X\ar[r]^{f_0}\ar[d] & X_{\K_0}\ar[d]\\
\spec(\K)\ar[r] & \spec\left(\K_0\right)
};
\end{displaymath}
moreover $f_0^{\ast}K_{X_{\K_0}}\cong K_X$, $K_{X_{\K_0}}$ is ample and $f_0^{\ast}\Omega^1_{X_{\K_0}}\cong\Omega^1_X$ (see \cite[Proposition 2.1.5]{AC:PhD}, \cite[\href{https://stacks.math.columbia.edu/tag/01V0}{tags 01V0} and \href{https://stacks.math.columbia.edu/tag/0D2P}{0D2P}]{TSP}). By \cite[Proposition 3]{L:SG}, $\Omega^1_X$ is $K_X$-semistable if and only if $\Omega^1_{X_{\K_0}}$ is $K_{X_{\K_0}}$-semistable. Accordingly one may assume $\K\subseteq\C$ without loss of generality. Consider the following Cartesian diagram
\begin{displaymath}
\xymatrix{
X_{\C}\ar[r]^{f}\ar[d] & X\ar[d]\\
\spec(\C)\ar[r] & \spec(\K)
}
\end{displaymath}
where $X_{\C}=X\times_{\spec(\K)}\spec(\C))$. Since $X_{\C}$ is a smooth, complex, projective variety (\cite[Lemma 2.1.3]{AC:PhD}) and $K_{X_{\C}}$ is ample, then $\Omega^1_{X_{\C}}$ is $K_{X_{\C}}$-semistable (see \cite[Theorem 1]{T:H}, for example). From all this, the claim follows.
\end{prf}
\begin{theorem}\label{th3.1}
The Higgs bundle $\fS$ is $K_X$-stable.
\end{theorem}
\begin{prf}
By \emph{Bertini's Theorem} (cfr. \cite[Corollary III.10.9 and Exercise III.11.3]{H:RC:1}), for $m\gg1$ there exist general, ample, smooth divisors $D_1,\dotsc,D_{n-1}\in\left|mK_X\right|$ such that $Y_p=D_1\cap\dotsc\cap D_p$ are smooth projective subvarieties of $X$ for any $p\in\{1,\dotsc,n-1\}$\footnote{
Moreover, these smooth projective varieties are all minimal and of general type. Indeed, by \emph{Adjunction Formula} $K_{Y_1}=\left(K_X+Y_1\right)_{\vert Y_1}$ (\cite[Exercise 21.5.B]{FOAG}) and this is an ample line bundle (\cite[Proposition 1.2.13 and Corollary 1.4.10]{L:RK}). Iterating this reasoning for all $Y_p$ one has the claim.
}; one puts $C=Y_{n-1}$. Since $\Omega^1_X$ is $K_X$-semistable, by \cite[Theorem 6.1]{M:R} $\Omega^1_{X\vert C}$ is semistable.\smallskip

\noindent Preliminary, one has the following inequality
\begin{displaymath}
0<\mu\left(S_{\vert C}\right)=\frac{1}{n+1}\int_Cc_1\left(S_{\vert C}\right)<\frac{1}{n}\int_Cc_1\left(\Omega^1_{X\vert C}\right)=\mu\left(\Omega^1_{X\vert C}\right).
\end{displaymath}
Let $\varpi\colon\fS_{\vert C}\to\fQ=\left(\cQ,\wvarphi\right)$ be a rank $s$ quotient Higgs bundle of $\fS$ and let $\cK=\ker(\varpi)$; $\cK$ is locally free because it is a kernel of an epimorphism of locally free sheaves on a smooth projective curve (see for example \cite[Proposition 1.1 and Corollary 1.4]{H:RC:2}). One sets
\begin{displaymath}
\cK_1=\cK\cap\Omega^1_{X\vert C},\cK_2=\cK\cap\cO_C.
\end{displaymath}
Consider the following commutative diagram
\begin{displaymath}
\xymatrix{
& 0\ar[d] & 0\ar[d] & 0\ar[d]\\
0\ar[r] & \cK_2\ar@{^{(}->}[d]\ar@{^{(}->}[r] & \cK\ar@{^{(}->}[d]\ar@{->>}[r] & \cK_1\ar@{^{(}->}[d]\ar[r] & 0\\
0\ar[r] & \cO_C\ar@{^{(}->}[r] & S_{\vert C}\ar@{->>}[r] & \Omega^1_{X\vert C}\ar[r] & 0%\\
% & & \cQ\ar[d]\ar@{->>}[r] & \cQ_0\ar[d]\ar[r] & 0\\
% & & 0 & 0
}
\end{displaymath}
where the rows are short exact sequences of sheaves. One has that $\cK_1$ is a sheaf of finite type on a Noetherian scheme, because it is a quotient of a locally free sheaf. In particular it is coherent, hence $\cK_2$ is a coherent sheaf as well, because it is a kernel of a morphism of coherent sheaves. On the other hand, $\cK_1$ is torsion-free, because it is a subsheaf of $\Omega^1_{X\vert C}$, hence $\cK_1$ is locally free. All this allows us to apply \cite[Proposition 1]{L:SG} thence $\cK_2$ is equal to either $\uzero_C$ or $\cO_C$.
\begin{enumerate}[1)]
\item Assume $\cK_2=\uzero_C$, then $s\in\{2,\dotsc,n\}$. By assumption $\varphi_{\vert C}(\cK)\subseteq\cK\otimes\Omega^1_C$ and by definition 
\begin{displaymath}
\varphi_{\vert C}(\cK)=\uzero_C\oplus\pi(\cK)\subseteq\left(\Omega^1_{X\vert C}\otimes\Omega^1_C\right)\oplus\Omega^1_C,
\end{displaymath}
hence $\cK$ is a Higgs subsheaf of $\fS_{\vert C}$ if and only if $\cK\subseteq\ker(\pi)\subseteq\Omega^1_{X\vert C}$. Let $Z$ be the subset of $\hGr_s\left(\fS_{\vert C}\right)$ parametrizing rank $s$ quotient Higgs bundles of $\fS_{\vert C}$ whose corresponding kernels are contained in $\ker(\pi)$. Let $U$ be an open subset of $C$ on which $S_{\vert U}$ is trivial, then $Z_{\vert U}$ is a subset of $\Gr_s\left(S_{\vert U}\right)=U\times_{\spec(\K)}\Gr_s\left(\K^{n+1}\right)$. More precisely, let $\left\{\ue_1,\dotsc,\ue_{n+1}\in\rH^0\left(S_{\vert C},U\right)\right\}$ be a local frame of $S_{\vert C}$ over $U$ such that $\ue_1$ generates $\cO_U$ and $\ue_3,\dotsc,\ue_{n+1}$ generate $\ker(\pi)_{\vert U}$. Let $I$ be an ordered subset of $\left\{2,\dotsc,n+1\right\}$ of cardinality $n-s+1$, then $Z_{\vert U}$ is contained in the union of ${U_I=U\times_{\spec(\K)}\left\{p_{\{1,\dotsc,n+1\}\setminus I}=0\right\}}$, where $\left\{p_{\alpha_1,\dotsc,\alpha_s}\mid 1\leq\alpha_1<\dotsc<\alpha_s\leq n+1\right\}$ are the \emph{Pl\"ucker coordinates} on $\Gr_s\left(\K^{n+1}\right)$. Then $Z_{\vert U}$ is a closed subscheme, because it is the intersection of $\hGr_s\left(\fS_{\vert U}\right)$ with the union of $U_I$'s; by gluing it follows that $Z$ is a closed subscheme. On the other hand, there exists a closed embedding of $Z_{red}$ in $\Gr_{s-1}\left(\Omega^1_{X\vert S}\right)$; indeed, since $\left\{p_{\{1,\dotsc,n+1\}\setminus I}=0\right\}$ is a \emph{Schubert cycle} corresponding to the $s$-dimensional quotients of $\K^{n+1}$ whose kernels do not contain $\ue_1$; by modding out these quotients by the subspace generated by $\ue_1$ we obtain $(s-1)$-dimensional quotients in $\K^n\cong\K^{n+1}/\langle\ue_1\rangle$, so that these $U_I$'s are isomorphic to ${U\times_{\spec(\K)}\Gr_{s-1}\left(\K^n\right)}$. As a consequence, the reduction of $Z$ is isomorphic to a closed subscheme of $\Gr_{s-1}\left(\Omega^1_{X\vert C}\right)$. Thus $\cQ^{\vee}$ corresponds to a subbundle of $T_{X\vert C}$, because $Z_{red}\hookrightarrow\Gr_{s-1}\left(\Omega^1_{X\vert C}\right)\cong\Gr_{n-s+1}\left(T_{X\vert C}\right)$ sends $[\cQ]$ to $\left[\cK^{\vee}\right]$. By the previous theorem
\begin{gather*}
\frac{1}{n}\deg\left(T_{X\vert C}\right)+\frac{1}{s}\deg(\cQ)\geq0\\
s\deg\left(T_{X\vert C}\right)+n\deg(\cQ)\geq0\\
s\deg\left(T_{X\vert C}\right)+(n+1)\deg(\cQ)>0
\end{gather*}
since $\deg(\cQ)>0$ and $\deg\left(T_{X\vert C}\right)<0$; in other words $\mu\left(S_{\vert C}\right)<\mu(\cQ)$.
\item Assume $\cK_2=\cO_C$. Under this hypothesis, it turns out that $\cQ$ is a quotient bundle of $\Omega^1_{X\vert C}$. Since $\Omega^1_{X\vert C}$ is semistable, $\mu\left(S_{\vert C}\right)<\mu\left(\Omega^1_{X\vert C}\right)\leq\mu(\cQ)$.
\end{enumerate}
From all this, $\fS$ is stable over the general  complete-intersection curve $C$ depending on $m$ and $K_X$, as explained at the begin of this proof, hence $\fS$ is $K_X$-stable.
\end{prf}
\noindent Recall
\begin{displaymath}
MY(X)=(-1)^n\int_X\left[c_2(X)-\frac{n}{2(n+1)}c_1(X)^2\right]\cdot K_X^{n-2};
\end{displaymath}
where $c_k(X)\stackrel{def.}{=}(-1)^kc_k\left(\Omega^1_X\right)$ for any $k\in\{1,\dotsc,n\}$. It is known that $MY(X)\geq0$ (cfr. \cite[Proposition 7.1]{M:Y:1}) which is called \emph{Miyaoka-Yau Inequality} (see \cite[Remark (i)]{Y:S-T}). By the previous theorem it follows another proof of this inequality.
\begin{theorem}[{cfr. \cite[Proposition 7.1]{M:Y:1}}]\label{th3.3}
The Miyaoka-Yau inequality $MY(X)\geq0$ holds.
\end{theorem}
\begin{prf}
It is enough to note that 
\begin{displaymath}
MY(X)=\int_X\Delta(S)\cdot c_1\left(K_X\right)^{n-2}
\end{displaymath}
and the claim follows from Theorem \ref{th1.1}.
\end{prf}

\subsection{Smooth projective varieties of general type such that $MY(X)=0$}

\noindent By Theorem \ref{th1.2} and Remark \ref{rem1.1}.\ref{rem1.1.a} the following corollary holds.
\begin{corollary}\label{cor3.1}
Let $X$ be a smooth projective variety of dimension $n\geq2$ with $K_X$ ample and $MY(X)=0$. The following statements hold:
\begin{enumerate}[a)]
\item\label{cor3.1.a} the Higgs bundle $\fS$ is curve semistable;
\item\label{cor3.1.b} $\displaystyle\Delta(S)=c_2(X)-\frac{n}{2(n+1)}c_1(X)^2\equiv_{num}0$;
\item\label{cor3.1.c} the Higgs bundle $\fS$ is semistable with respect to any polarization of $X$.
\end{enumerate}
\end{corollary}
\noindent Now we are in position to study the positivity conditions for $\fS$.
\begin{lemma}\label{lem3.1}
Let $X$ be a smooth projective variety of dimension $n\geq2$ with $K_X$ ample and $MY(X)=0$. Then $\fS$ is a weakly H-ample Higgs bundle over $X$.
\end{lemma}
\begin{prf}
Since $\fS$ is curve semistable (Corollary \ref{cor3.1}.\ref{cor3.1.a})
\begin{displaymath}
\mu_{\min}\left(f^{\ast}\fS\right)=\mu\left(f^{\ast}\fS\right)=\frac{1}{n+1}\int_Cf^{\ast}c_1(S)=\delta\int_Cf^{\ast}c_1\left(K_X\right)
\end{displaymath}
for all $f\colon C\to X$. By the \emph{Weak H-Ampleness Criterion} (Theorem \ref{th2.1}), taking $h=K_X$ and $\displaystyle\delta=\frac{1}{n+1}$, it follows the claim.
\end{prf}
\begin{remark}
$S$, as an ordinary vector bundle, is neither stable nor ample. Indeed, $\cO_X$ is a quotient bundle of $S$ and $\mu(S)>\mu(\cO_X)=0$.
\end{remark}
\noindent By the proof of Theorem \ref{th3.1}, by applying Proposition \ref{prop2.1}.\ref{prop2.1.d}, one has the following corollary.
\begin{corollary}\label{cor3.2}
Let $X$ be a smooth projective variety of dimension $n\geq2$ with $K_X$ ample and $MY(X)=0$. Then the cotangent bundle $\Omega^1_X$ is ample.
\end{corollary}
\begin{prf}
By the premises and Remark \ref{rem2.1}.\ref{rem2.1.a}, it is enough to prove that $\psi$, the Higgs field induced by $\varphi$ on the quotient Higgs bundle $\Omega^1_X$, vanishes. Reasoning on the stalks, one has
\begin{displaymath}
\forall x\in X,\widetilde{\omega}\in\Omega^1_{X,x},s\in\cO_{X,x},\varphi_x\left(\widetilde{\omega},s\right)=\left(0,\widetilde{\omega}\right)\in\left(\Omega^1_{X,x}\right)^{\otimes2}\oplus\Omega^1_{X,x}.
\end{displaymath}
Since $\psi_x\circ pr_{1,x}=qr_{1,x}\circ\varphi_x=0$, where $pr_{1,x}\colon\Omega^1_{X,x}\oplus\cO_{X,x}\to\Omega^1_{X,x}$ and ${qr_{1,x}\colon\left(\Omega^1_{X,x}\right)^{\otimes2}\oplus\Omega^1_{X,x}\to\left(\Omega^1_{X,x}\right)^{\otimes2}}$ are the canonical projections, then $\psi_x$ vanishes for any $x\in X$, \emph{i.e.} $\psi$ vanishes.
\end{prf}
\noindent Moreover, the following proposition holds.
\begin{proposition}\label{prop3.1}
Let $X$ be a smooth projective variety of dimension $n\geq3$ with $K_X$ ample and $MY(X)=0$. Then $\fS$ is not strongly H-ample.
\end{proposition}
\begin{prf}
Since $\Omega_X^1$ is a rank $n$ ample quotient bundle of $\fS$, by Remark \ref{rem2.3} and Proposition \ref{prop2.1}.\ref{prop2.1.d} one has the claim.
\end{prf}

\section{On twisted Simpson system}

\subsection{Premise}

\noindent As it is customary, we consider twistings by rational divisors formally, these make sense after pulling back to a (possibly ramified) finite covering of $X$. To be clear, we recall the following theorem.
\begin{theorem}[{\cite[Lemma 2.1]{B:G}}, cfr. also {\cite[Theorem 4.1.10]{L:RK}}]
Let $X$ be a projective variety, let $L$ be a line bundle over X, and fix any positive integer number $d$. Then there exist a projective variety $Y$, a finite surjective map $f\colon Y\to X$, together with a line bundle $M$ over $Y$ such that $f^{\ast}L=M^{\otimes d}$.
\end{theorem}
\noindent Let $L\in\Pic(X)\otimes_{\Z}\Q$. There exists $d\gg1$ such that $L^{\otimes d}\in\Pic(X)$ so one defines
\begin{displaymath}
c_1(L)=\frac{1}{d}c_1\left(L^{\otimes d}\right)\in\rA^1(X)_{\Q};
\end{displaymath}
in this way, one extends by $\Q$-linearity the results of \cite[Section 2.5]{F:W} and it makes sense to consider the degree of $\Q$-line bundles over $X$ with respect to a polarization of $X$. On the other hand, as it customs, one says that $L$ is ample (nef, respectively) over $X$ if $M$ is an ample (a nef, respectively) line bundle over $Y$, where $f^{\ast}L^{\otimes d}=M^{\otimes d}$; such $f$, $M$ and $Y$ exist by the previous theorem of course. Moreover, it follows that the properties of being weakly H-ample and H-nef are invariant under such coverings (see Propositions \ref{prop2.1}.\ref{prop2.1.c} and \ref{prop2.2}.\ref{prop2.2.a}).

\subsection{Main results}

\noindent More generally than Lemma \ref{lem3.1}, the following theorem holds.
\begin{theorem}\label{th4.1}
Let $X$ be a smooth projective variety of dimension $n\geq2$ with $K_X$ ample and $MY(X)=0$. Then the Higgs bundle $\fS_{\beta}=\fS\left(-\beta K_X\right)$ is:
\begin{enumerate}[a)]
\item\label{th4.1.a} H-nef for every rational number $\displaystyle\beta\leq\frac{1}{n+1}$;
\item\label{th4.1.b} weakly H-ample for every rational number $\displaystyle\beta<\frac{1}{n+1}$.
\end{enumerate}
\end{theorem}
\begin{prf}
By Theorem \ref{th3.1}, after twisting by a $\Q$-line bundle, mimicking the proof \cite[Lemma 6.4.4]{L:RK} one has that $\fS_{\beta}$ is curve semistable, \emph{i.e.} for every morphism $f\colon C\to X$, the pullback Higgs bundle $f^{\ast}\fS_{\beta}$ is semistable, hence
\begin{gather*}
\mu\left(f^{\ast}\fS_{\beta}\right)=\mu_{\min}\left(f^{\ast}\fS_{\beta}\right)=\frac{1}{n+1}\int_Cf^{\ast}c_1\left(S_{\beta}\right)=\frac{1}{n+1}\int_{\olC}c_1\left(S\left(-\beta K_X\right)_{\vert\olC}\right)\\
=\frac{1}{n+1}\int_{\olC}\left(1-\beta(n+1)\right)c_1\left(K_{X\vert\olC}\right)=\left(\frac{1}{n+1}-\beta\right)\int_{\olC}K_{X\vert\olC}\geq0
\end{gather*}
where $\olC=f(C)$; the last inequality holds as $K_X$ is ample. By the properties of the HN-filtration of Higgs bundles, any quotient Higgs bundle of $f^{\ast}\fS_{\beta}$ has non-negative degree. Furthermore, by the previous computation
\begin{displaymath}
\int_Cf^{\ast}c_1\left(S_{\beta}\right)=(1-(n+1)\beta)\int_{\olC}K_{X\vert\olC}\geq0.
\end{displaymath}
Thus, by the \emph{weak H-ampleness} and \emph{H-nefness Criteria} (Theorems \ref{th2.1} and \ref{th2.2}, respectively), one has the ampleness and nefness of $\det\left(\fS_{\beta}\right)$ as a function of $\beta$, so that the claim follows.
\end{prf}
\begin{remark}
For any rational number $\displaystyle0<\beta\leq\frac{1}{n+1}$, $\fS_{\beta}$ is a H-nef Higgs bundle whose underlying vector bundle $S_{\beta}$ is not nef; indeed, $-\beta K_X$ is a quotient $\Q$-line bundle of $S_{\beta}$ which is not nef.
\end{remark}
\noindent In consequence, we generalize \cite[Corollary at page 294]{M:Y:3} to higher dimensional case.
\begin{corollary}\label{cor4.1}
Let $X$ be a smooth projective variety of dimension $n\geq2$ with $K_X$ ample and $MY(X)=0$. Then the vector bundle $\Omega_{\beta}=\Omega^1_X\left(-\beta K_X\right)$ is:
\begin{enumerate}[a)]
\item nef for every rational number $\displaystyle\beta\leq\frac{1}{n+1}$;
\item ample for every rational number $\displaystyle\beta<\frac{1}{n+1}$.
\end{enumerate}
Consequently for each irreducible projective curve $C$ on $X$, there exists an effective divisor $D$ on the normalization $\wtC$ of $C$ such that the inequality $\displaystyle2p_a(C)-2\geq\frac{1}{n+1}\int_CK_{X\vert C}+\deg(D)$ holds, where $p_a(C)$ is the \emph{arithmetic genus of $C$}. In particular $X$ does not contain irreducible curves of arithmetic genus neither $0$ nor $1$.
\end{corollary}
\begin{prf}
$\left(\Omega_{\beta},0\right)$ is a quotient Higgs bundle of $\fS_{\beta}$ hence it is H-nef (Proposition \ref{prop2.2}.\ref{prop2.2.b} and Corollary \ref{cor3.2}) and then nef in the usual sense (Remark \ref{rem2.1}.\ref{rem2.1.a}). Let $C\subsetneqq X$ be an irreducible projective curve, let $\left(\wtC,\nu\right)$ be its normalization, and let $\iota\colon\wtC\xrightarrow{\nu}C\hookrightarrow X$. One has a morphism
\begin{displaymath}
\xymatrix{
\iota^{\ast}\Omega_{\textstyle\frac{1}{n+1}}\ar[r] & \displaystyle\Omega^1_{\wtC}\left(-\frac{1}{n+1}\iota^{\ast}K_X\right)=\cO_{\wtC}\left(K_{\wtC}-\frac{1}{n+1}\iota^{\ast}K_X\right)
}
\end{displaymath}
which does not need to be surjective at ramification points. So one considers the saturation of the image of this morphism: it is a line bundle $\displaystyle\cO_{\wtC}\left(K_{\wtC}-\frac{1}{n+1}\iota^{\ast}K_X-D\right)$ for some effective divisor $D$ on $\wtC$. Since this is a quotient bundle of $\displaystyle\iota^{\ast}\Omega_{\textstyle\frac{1}{n+1}}$, it is nef. By all this
\begin{displaymath}
2p_a(C)-2\geq2p_a\left(\wtC\right)-2=\deg\Omega^1_{\wtC}\geq\frac{1}{n+1}\int_CK_{X\vert C}+\deg(D)>0.
\end{displaymath}
As a consequence $X$ contains irreducible curves of arithmetic genus neither $0$ nor $1$.
\end{prf}
\noindent On the other hand, the previous inequality is not sharp as the following example proves.
\begin{example}
Let $X$ be a \emph{fake complex projective plane}, \emph{i.e.} $X$ is a smooth projective surface with the same \emph{Betti numbers of} $\bP^2_{\C}$ but is not isomorphic to $\bP^2_{\C}$. Its canonical bundle $K_X$ is ample hence $X$ is a minimal smooth surface of general type, $\displaystyle{\int_X3c_2(X)=\int_Xc_1(X)^2=9}$ and the \emph{Picard number of $X$ is} $1$. Let $\cO_X(1)$ be the ample generator of the torsion-free part of $\Pic(X)$ such that $K_X=\cO_X(3)$ and $c_1\left(\cO_X(1)\right)^2=1$ (cfr. \cite[Section 1.1]{DiC:LF}). Moreover, let $f\colon C\to X$ be non-constant then $g(C)>2$ by \cite[Lemma 2.2]{K:JH} and \cite[Proposition 2.3]{DiC:LF}\footnote{
To be clear, whether $C$ is contained in $X$ or is the normalization of a reduced irreducible curve in $X$, the genus of $C$ is greater than $2$.
}. However, for such a curve $C$ the previous corollary predicts:
\begin{displaymath}
2g(C)-2\geq\frac{1}{2+1}c_1\left(\cO_X(3)\right)\cdot C=c_1\left(\cO_X(1)\right)\cdot C>0\iff g(C)>1.
\end{displaymath}
\end{example}
\noindent A little recap: in Lemma \ref{lem3.1}, Theorem \ref{th4.1}, Corollaries \ref{cor3.1} and \ref{cor4.1} $X$ is a smooth projective variety of dimension $n\geq2$ with $K_X$ ample and $MY(X)=0$; under these hypotheses it follows that:
\begin{enumerate}[a)]
\item $\fS$ is curve semistable;
\item $\fS_{\beta}$ is H-nef hence $\Omega^1_{\beta}$ is nef for every rational number $\displaystyle\beta\leq\frac{1}{n+1}$;
\item $\fS_{\beta}$ is weakly H-ample hence $\Omega^1_{\beta}$ is ample for every rational number $\displaystyle\beta<\frac{1}{n+1}$;
\item for each irreducible projective curve $C$ on $X$, there exists an effective divisor $D$ on the normalization $\wtC$ of $C$ such that the inequality $\displaystyle2p_a(C)-2\geq\frac{1}{n+1}\int_CK_{X\vert C}+\deg(D)$ holds;
\item $X$ contains irreducible curves of arithmetic genus neither $0$ nor $1$.
\end{enumerate}
These results can be generalised, partially.
\begin{lemma}\label{lem4.1}
Let $X$ be a smooth projective variety of general type of dimension $n\geq2$, \emph{i.e.} $K_X$ is big and nef. If the Higgs bundle $\fS$ is semistable with respect to some polarization of $X$ and $\Delta(S)\equiv_{num}0$ then:
\begin{enumerate}[a)]
\item $\fS$ is curve semistable;
\item\label{lem4.1.b} $\fS_{\beta}$ is H-nef hence $\Omega^1_{\beta}$ is nef for every rational number $\displaystyle\beta\leq\frac{1}{n+1}$;
\item\label{lem4.1.c} for each irreducible projective curve $C$ on $X$, there exists an effective divisor $D$ on the normalization $\wtC$ of $C$ such that the inequality $\displaystyle2p_a(C)-2\geq\frac{1}{n+1}\int_CK_{X\vert C}+\deg(D)$ holds;
\item\label{lem4.1.d} $X$ contains irreducible curves of arithmetic genus neither $0$ nor $1$;
\item $K_X$ is ample;
\item\label{lem4.1.f} $\fS_{\beta}$ is weakly H-ample hence $\Omega^1_{\beta}$ is ample for every rational number $\displaystyle\beta<\frac{1}{n+1}$.
\end{enumerate}
\end{lemma}
\begin{prf}
$\fS$ is curve semistable by Theorem \ref{th1.2}.\ref{th1.2.a} and by applying Proposition \ref{prop2.2}.\ref{prop2.2.b}. Repeating the proof of Theorem \ref{th4.1}.\ref{th4.1.a}, one has the statement \ref{lem4.1.b}. Repeating the proof of Corollary \ref{cor4.1}, one has the statement \ref{lem4.1.c}; in particular, $X$ does not contain rational curves, hence $K_X$ is ample by \cite[Exercise 7.13.8]{D:O}. Thus and by statement \ref{lem4.1.c} one completes the proof of statement \ref{lem4.1.d}. Repeating the proof of Theorem \ref{th4.1}.\ref{th4.1.b} and applying Proposition \ref{prop2.1}.\ref{prop2.1.d}, one has the statement \ref{lem4.1.f}.
\end{prf}

\section{On the algebraic hyperbolicity of smooth projective varieties\\ with ample canonical bundle}

\noindent For clarity, we recall the following definition.
\begin{definition}[cfr. {\cite[Definition 6.3.23.(i)]{L:RK}}]\label{def5.1}
A smooth projective variety $Y$ is \emph{algebraically hyperbolic} if there exists an ample divisor $L$ on $Y$ and $\epsilon\in\R_{>0}$ such that the arithmetic genus $p_a(C)$ of any projective curve $C$ on $Y$ satisfies the inequality
\begin{displaymath}
2p_a(C)-2\geq\epsilon\int_CL_{\vert C}.
\end{displaymath}
\end{definition}
\noindent By Corollary \ref{cor4.1}, using the notations of the previous definition, the following theorem holds.
\begin{theorem}\label{th3.5}
Let $X$ be a smooth projective variety of dimension $n\geq2$ with $K_X$ ample and $MY(X)=0$. Then $X$ is algebraically hyperbolic choosing $L=K_X$ and $\displaystyle\epsilon\geq\frac{1}{n+1}$.
\end{theorem}
\begin{remark}
The previous theorem follows also by Corollary \ref{cor3.2} and \cite[Theorem 6.3.26]{L:RK}. But in this way one does not obtain any estimate of $\epsilon$.
\end{remark}
\noindent Finally, by \cite[Proposition 9.8]{S:CT:1} one has another proof of the following theorem.
\begin{theorem}[{cfr. \cite[Theorem 4 and Remark (iii)]{Y:S-T}}]\label{th5.1}
Let $X$ be a smooth projective variety of dimension $n\geq2$ with $K_X$ ample and $MY(X)=0$. If $\K=\C$ then $X$ is uniformized by $\B^n$ (the unit ball of $\C^n$).
\end{theorem}
\begin{remark}
By \cite[Example 2.1.2.2]{Z:L}, if $\K=\C$ then Corollary \ref{cor3.2} can be viewed as a consequence of the previous Theorem.
\end{remark}
\noindent In conclusion, I posit the following conjecture.
\begin{conjecture}\label{conj}
Let $X$ be a smooth projective variety of dimension $n\geq2$ with $K_X$ ample, let $\displaystyle\alpha=\frac{\int_Xc_2(X)\cdot K_X^{n-2}}{\int_Xc_1^2(X)\cdot K_X^{n-2}}\in\Q_{>0}$. Then $X$ is algebraically hyperbolic with $L=K_X$ and $\displaystyle\epsilon\geq1-2\alpha$.
\end{conjecture}
\begin{remark}
If $n=2$ then this conjecture holds. Indeed, $\Omega_{\beta}$ is nef for any rational number $\beta\leq 1-2\alpha$ by \cite[Remark at page 294]{M:Y:1}; one repeats the proof of Corollary \ref{cor4.1} and has the claim accordingly.
\end{remark}
\noindent{\bf Acknowledgment.} I am deeply grateful to my Ph.D. advisors, Ugo Bruzzo and Beatriz Gra\~na Otero, for their invaluable guidance, energy and unwavering support. I would also like to thank SISSA for its hospitality while portion of this work was being completed. I am very grateful to Angelo Felice Lopez for pointing out several oversights in the references, and to Indranil Biswas, Snehajit Misra and Nabanita Ray for bringing their paper to my attention. Finally, special thanks go to the anonymous referee for their insightful suggestions.\medskip

\noindent {\bf AI declaration.} No artificial intelligence software was used in the writing of this article.\medskip

\noindent {\bf Statement about competing or financial interests.} The author has no competing or financial interests to declare that are relevant to the content of this article.

\end{document}